\documentclass[a4paper,12pt]{article}
\usepackage{amssymb}
\usepackage{amsmath}
\usepackage{graphicx}
\newcommand{\C}{\mathbb{C}}
\newcommand{\coleq}{\mathrel{:=}}
\newcommand{\der}{\mathrm{d}}
\newcommand{\Mn}{\mathbb{E}}
\newcommand{\N}{\mathbb{N}}
\newcommand{\nep}{\mathrm{e}}
\newcommand{\Pb}{\mathbb{P}}
\newcommand{\Rl}{\mathbb{R}}
\newcommand{\tit}[1]{\noindent\textbf{#1}}
\newcommand{\un}{\mathbf{1}}

\newtheorem{thm}{Theorem}
\newtheorem{lem}[thm]{Lemma}
\newtheorem{defi}[thm]{Definition}
\newtheorem{rmk}[thm]{Remark}

\newtheorem{xpl}[thm]{Example}
\newtheorem{prop}[thm]{Proposition}
\title{A concentration theorem for the equilibrium measure of Markov chains with nonnegative coarse Ricci curvature}
\author{Laurent Veysseire}
\date{}
\begin{document}
\maketitle
\abstract{In this article, we prove a concentration inequality of the order of the exponential of a double integral of the coarse Ricci curvature for the equilibrium measure of a Markov chain, in the case when this curvature is nonnegative.
This is, to the author's knowledge, the first concentration result in a discrete setting using a non-constant curvature instead of its infimum.}

\section*{Introduction}

For a Markov chain on a Polish space, a nonnegative coarse Ricci curvature means that the distributions after one step of the chain are closer (in the sense of the $W_1$ distance) than their starting points are \cite{OLLI09}.
Remind that the $W_1$ (Wasserstein) metric between two probability measures is the infimum over the set of couplings between this two probability measures of the expectation of the distance between the two points.

In the case when the space is $\varepsilon$-geodesic (see Definition \ref{gd}), a nonnegative coarse Ricci curvature allows to extend the local attractiveness of a point $x_0$ to a global one (see \cite{OLLI09} or Lemma \ref{lem1}).
The attractiveness of a point implies exponential concentration of the equilibrium probability measure around this point, if the Markov chain does not spread out too quickly.

One of the simplest example is the random walk on $\N$ where we jump from $n$ to $n+1$ with probability $p$ and to $(n-1)_+$ with probability $1-p$.
In this case, the coarse Ricci curvature is $0$.
If $p<\frac{1}{2}$, then $0$ is attractive and we have exponential concentration.
If $p\geq\frac{1}{2}$, then we don't have any attractive point, neither do we have any invariant probability measure.

Here we prove that the concentration of the equilibrium measure around an attractive point behaves at least like the exponential of a double integral of the coarse Ricci curvature.

We may remark that this is the right behaviour of the invariant distribution for diffusion processes on the real line, as we see in the example below.

\begin{xpl}Let us consider a diffusion process on the real line whose generator takes the form:
$$Lf=\frac{\der^2f}{\der x^2}-\frac{\der V}{\der x}\frac{\der f}{\der x}$$
where the energy $V(x)$ is smooth.
Then the coarse Ricci curvature is $\frac{\der^2V}{\der x^2}$, and the measure $\nep^{-V(x)}\der x$ is reversible.
We see that the density of the invariant measure is exacly a double integral of the coarse Ricci curvature.
\end{xpl}

\section{The concentration Theorems}

We define $\varepsilon$-geodesic spaces as in \cite{OLLI09}.

\begin{defi}\label{gd} Let $\varepsilon>0$.
A metric space $(X,d)$ is said to be $\varepsilon$-geodesic if for each $(x,y)\in X^2$, there exists $n\in\N$ and a sequence $x=x_0,x_1,\ldots,x_n=y\in X$ such that $d(x_i,x_{i+1})\leq\varepsilon$ for each $0\leq i\leq n-1$ and $d(x,y)=\sum_{i=0}^nd(x_i,x_{i+1})$.
\end{defi}
For a Markov chain with transition kernel $P$ on a $\varepsilon$-geodesic space, we will denote by $K_\varepsilon(x)$ the local coarse Ricci curvature at $x$: $$K_{\varepsilon}(x)\coleq\inf_{y\in X|0<d(x,y)\leq\varepsilon}\kappa(x,y)$$
Where $\kappa(x,y)\coleq 1-\frac{W_1(P_x,P_y)}{d(x,y)}$ is the coarse Ricci curvature between $x$ and $y$ as defined in \cite{OLLI09}.

Here we will prove the following concentration result for the equilibrium measure of Markov Chains:

\begin{thm}\label{tm1}Let $X$ be an $\varepsilon$-geodesic metric space and $P$ be the transition kernel of a Markov chain on $X$.
Assume that:
\begin{itemize}
\item there exists $\rho>0$ and a point $x_0$ such that $x_0$ is $\rho$-attractive for the Markov chain in the sense that
$$\forall x|\varepsilon<d(x,x_0)\leq2\varepsilon,W_1(\delta_{x_0},P_x)\leq d(x_0,x)-\rho,$$
\item  there exists a non-increasing function $K:\Rl_+\mapsto\Rl_+$ satisfying:
$$K_\varepsilon(x)\geq K(d(x,x_0)),$$
\item there exists $s>0$ such that for any $x\in X$, any $1$-lipschitz function $f:X\mapsto\Rl$ and any $\lambda\in\Rl$, we have:
$$\Mn_{P_x}\left[\nep^{\lambda f}\right]\leq \nep^{\lambda\Mn_{P_x}[f]+\frac{\lambda^2s^2}{2}}.$$
\end{itemize}
Then we have, for every $l>2\varepsilon+\frac{\ln(2)s^2}{\rho}$ and any equilibrium measure $\pi$:
$$\Pb_{x\sim\pi}(d(x,x_0)\geq l)\leq C_0\nep^{-\frac{1}{2s^2}\Phi(l)}$$
with
$$\Phi(l)\coleq\rho l+\int_{2\varepsilon}^l\left(\int_{2\varepsilon}^uK(v)\der v\right)\der u$$
and
$$C_0=\frac{\nep^{\frac{3\varepsilon}{2s^2}\max(3\varepsilon,\rho+\frac{\ln(2)s^2}{\rho})-\frac{\rho^2}{4s^2}+\frac{1}{2s^2}\left(\rho(2\varepsilon+\frac{\ln(2)s^2}{\rho})+\int_{2\varepsilon}^{2\varepsilon+\frac{\ln(2)s^2}{\rho}}\int_{2\varepsilon}^uK(v)\der v\der u\right)}}{1-\nep^{-\frac{\rho^2}{4s^2}}}$$
\end{thm}

\begin{rmk} If $K=0$, we obtain exponential concentration, as proved in \cite{OLLI09}.
\end{rmk}

\begin{prop}\label{tm2} If closed balls are compact, then under the hypotheses of Theorem \ref{tm1}, there exists an equilibrium measure.
\end{prop}

\begin{rmk} In the case when $\int_0^\infty K(r)\der r=\infty$, and for some (hence any) $x_0\in X$, $W_1(\delta_{x_0},P_{x_0})<\infty$, then for any $\rho>0$, there exists a $\varepsilon>0$ large enough such that $x_0$ is $\rho$-attractive.
This is a trivial consequence of the Lemma below.
\end{rmk}

\begin{lem}\label{lem1} Let $X$ be an $\varepsilon$-geodesic metric space and $P$ be the transition kernel of a Markov chain such that there exists a non-increasing function $K:\Rl_+\mapsto\Rl_+$ and a point $x_0\in X$ satisfying:
$$K_\varepsilon(x)\geq K(d(x,x_0)).$$
Then we have $$\Mn_{y\sim P(x)}[d(x_0,y)]\leq d(x_0,x)-F(d(x_0,x)),$$ where
$$F(l)\coleq\left\{\begin{array}{ll}\rho+\int_{2\varepsilon}^lK(u)\der u & \textrm{if }2\varepsilon\leq l\\
\rho & \textrm{if }\varepsilon<l\leq 2\varepsilon\\
-J(0) & \textrm{if }l\leq\varepsilon\end{array}\right.$$
with $\rho\coleq\inf_{x|\varepsilon<d(x,x_0)\leq2\varepsilon} d(x,x_0)-W_1(P(x),\delta_{x_0})$ and $J(x_0)=W_1(P(x_0),\delta_{x_0})$.
\end{lem}

\tit{Proof :}

If $\varepsilon< d(x,x_0)\leq 2\varepsilon$, this is just the definition of $\rho$.
If $d(x,x_0)\leq\varepsilon$, we have $K_\varepsilon(x_0)\geq0$, so $W_1(P(x),\delta_{x_0})\leq W_1(P(x),P(x_0))+W_1(P(x_0),\delta_{x_0})\leq d(x,x_0)+J(x_0)$.
If $d(x,x_0)\geq 2\varepsilon$, there exists $x_1,\ldots,x_n=x$ such that $d(x_i,x_{i+1})\leq\varepsilon$, $\varepsilon<d(x_1,x_0)\leq2\varepsilon$ and $d(x,x_0)=d(x_1,x_0)+\sum_{i=1}^{n-1}d(x_i,x_{i+1})$.
We have then 
\begin{align*}W_1(P(x),\delta_{x_0})&\leq W_1(P(x_1),\delta_{x_0})+\sum_{i=1}^{n-1}W_1(P(x_i),P(x_{i+1}))\\
{}&\leq d(x_1,x_0)-\rho+\sum_{i=1}^{n-1}(1-K(d(x_i,x_0)))d(x_i,x_{i+1})\\
{}&\leq d(x,x_0)-\rho-\sum_{i=1}^{n-1}\int_{d(x_i,x_0)}^{d(x_{i+1},x_0)}K(l)\der l\\
{}&\leq d(x,x_0)-\rho-\int_{d(x_1,x_0)}^{d(x,x_0)}K(l)\der l\\
{}&\leq d(x,x_0)-F(d(x,x_0)).\square\end{align*}

\begin{lem}\label{lem2}
Let $\mu$ be a probability measure on $X$ and $s>0$ be such that for any $1$-lipschitz function $f$, we have the following inequality:
$$\Mn_\mu\left[\nep^{\lambda f}\right]\leq \nep^{\lambda\Mn_\mu[f]+\frac{\lambda^2s^2}{2}}.$$
Then, for each $\mathcal{C}^1$ function $g:\Rl\mapsto\Rl$ such that $g'$ is Lipschitz and $\|g'\|_{lip}<\frac{1}{s^2}$ and for each $1$-lipschitz function $f$, we have:
$$\Mn_\mu\left[\nep^{g\circ f}\right]\leq \frac{\nep^{g\left(\Mn_\mu[f]\right)+\frac{s^2g'^2\left(\Mn_\mu[f]\right)}{2\left(1-s^2\|g'\|_{lip}\right)}}}{\sqrt{1-s^2\|g'\|_{lip}}}.$$
\end{lem}

\tit{Proof :}

For each $x\in X$, we have
$$\nep^{g\circ f(x)}\leq \nep^{g\left(\Mn_\mu[f]\right)+\left(f(x)-\Mn_\mu[f]\right)g'\left(\Mn_\mu[f]\right)+\frac{\left(f(x)-\Mn_\mu[f]\right)^2}{2}\|g'\|_{lip}}.$$
Now we use the fact that the Laplace transform of a Gaussian measure $\mathcal{N}(M,\sigma^2)$ is:
$$\int_{-\infty}^{\infty}\nep^{\lambda u}\nep^{-\frac{(u-M)^2}{2\sigma^2}}\frac{\der u}{\sqrt{2\pi\sigma^2}}=\nep^{\lambda M+\frac{\lambda^2\sigma^2}{2}}.$$
So, taking $\lambda=f(x)-\Mn_\mu[f]$, $M=g'\left(\Mn_\mu[f]\right)$ and $\sigma^2=\|g'\|_{lip}$, we get:
$$\nep^{g(f(x))}\leq\nep^{g\left(\Mn_\mu[f]\right)}\int_{-\infty}^{\infty}\nep^{u\left(f(x)-\Mn_\mu[f]\right)}\nep^{-\frac{(u-g'\left(\Mn_\mu[f]\right))^2}{2\|g'\|_{lip}}}\frac{\der u}{\sqrt{2\pi\|g'\|_{lip}}}.$$
Integrating this inequality with respect to $\mu$ and using our assumption yields:
$$\Mn_\mu\left[\nep^{g\circ f}\right]\leq\nep^{g\left(\Mn_\mu[f]\right)}\int_{-\infty}^{\infty}\nep^{\frac{u^2s^2}{2}}\nep^{-\frac{(u-g'\left(\Mn_\mu[f]\right))^2}{2\|g'\|_{lip}}}\frac{\der u}{\sqrt{2\pi\|g'\|_{lip}}}=\nep^{g\left(\Mn_\mu[f]\right)}\frac{\nep^{\frac{s^2g'^2\left(\Mn_\mu[f]\right)}{2(1-s^2\|g'\|_{lip}}}}{\sqrt{1-s^2\|g'\|_{lip}}}$$
as needed.$\square$

\begin{thm}\label{princ} Let $X$ be a $\varepsilon$-geodesic metric space and $P$ be the transition kernel of a Markov chain.
Assume that there exists a non-increasing function $K:\Rl_+\mapsto\Rl_+$ and a point $x_0\in X$ satisfying:
$$K_\varepsilon(x)\geq K(d(x,x_0))$$
and that there exists $s>0$ such that for any $x\in X$, any $1$-lipschitz function $f:X\mapsto\Rl$ and any $\lambda\in\Rl$, we have:
$$\Mn_{P_x}\left[\nep^{\lambda f}\right]\leq \nep^{\lambda\Mn_{P_x}[f]+\frac{\lambda^2s^2}{2}}.$$
Let $F$ be defined as in Lemma \ref{lem1}.
Then, for every pair $(\alpha,d_0)\in\Rl_+^2$ satisfying:
\begin{itemize}
\item $d_0\geq2\varepsilon$
\item $F(d_0)>\frac{s^2K(d_0)}{2}$
\item $\alpha<\frac{1}{s^2K(d_0)}$
\item $\displaystyle C_{\alpha,d_0}\coleq \frac{\nep^{-\alpha F(d_0)^2\left(1-\frac{\alpha s^2}{2(1-\alpha s^2K(d_0))}\right)}}{\sqrt{1-\alpha s^2K(d_0)}}<1$
\end{itemize}
we have the following concentration inequality for any equilibrium measure $\pi$ of the Markov chain and any $l\geq d_0$:
$$\Pb_{x\sim\pi}(d(x,x_0)\geq l)\leq C'_{\alpha,d_0}\frac{C_{\alpha,d_0}}{1-C_{\alpha,d_0}}\nep^{-\alpha(\varphi(l)-\varphi(d_0))}$$
where $\varphi(l)=\int_0^lF(u)\der u$, and $C'_{\alpha,d_0}\coleq\nep^{\alpha\un_{J(x_0)+\varepsilon>d_0-F(d_0)}\int_{d_0-F(d_0)}^{J(x_0)+\varepsilon}\sup(F(d_0),F(u))\der u}$.
\end{thm}

\tit{Proof :} we set $\psi(x)=\alpha\varphi(x)$ if $x\geq d_0$ and $\psi(x)=\alpha(\psi(d_0)-(d_0-x)F(d_0))$ if $x<d_0$.
Under our assumptions, $\psi$ is convex and increasing, and we have $\|\psi'\|_{lip}=\alpha K(d_0)<\frac{1}{s^2}$.
Our goal is to bound the quantity $\Mn_{x\sim\pi}\left[\nep^{\psi(d(x,x_0))}\un_{d(x,x_0)\geq d_0}\right]$.
We have:
\begin{align*}\Mn_{x\sim\pi}\left[\nep^{\psi(d(x,x_0))}\un_{d(x,x_0)\geq d_0}\right]&=\Mn_{x\sim\pi}\left[\Mn_{y\sim P_x}\left[\nep^{\psi(d(y,x_0))}\un_{d(y,x_0)\geq d_0}\right]\right]\\
&\leq\Mn_{x\sim\pi}\left[\Mn_{y\sim P_x}\left[\nep^{\psi(d(y,x_0))}\right]\right].\end{align*}
Using Lemma \ref{lem2} with $\mu=P_x$ and $g=\psi$, and Lemma \ref{lem1}, we get:
\begin{multline*}\Mn_{x\sim\pi}\left[\nep^{\psi(d(x,x_0))}\un_{d(x,x_0)\geq d_0}\right]\leq\\
\Mn_{x\sim\pi}\left[\frac{\nep^{\psi(d(x,x_0)-F(d(x,x_0)))+\frac{\alpha^2s^2F(d_0)^2}{2(1-\alpha s^2K(d_0))}}}{\sqrt{1-\alpha s^2K(d_0)}}\un_{d(x,x_0)<d_0}\right.\\
\left.+\frac{\nep^{\psi(d(x,x_0)-F(d(x,x_0)))+\frac{\alpha^2s^2F(d(x,x_0))^2}{2(1-\alpha s^2K(d_0))}}}{\sqrt{1-\alpha s^2K(d_0)}}\un_{d(x,x_0)\geq d_0}\right].\end{multline*}
The function $l\mapsto l-F(l)$ is nondecreasing on $[0,\varepsilon]$ and on $(\varepsilon,d_0)$, and $\psi$ is an increasing function.
Then, for $d(x,x_0)<d_0$, we have $\psi(d(x,x_0)-F(d_0))\leq\psi(\max(J(x_0)+\varepsilon,d_0-F(d_0))=\ln(C'_{\alpha,d_0})+\alpha(\varphi(d_0)-F^2(d_0))$.

For $d(x,x_0)\geq d_0$, we have $\psi(d(x,x_0)-F(d(x,x_0)))\leq \psi(d(x,x_0))-\alpha F^2(d_0)$.
So we get:
$$\Mn_{x\sim\pi}\left[\nep^{\psi(d(x,x_0))}\un_{d(x,x_0)\geq d_0}\right]\leq C'_{\alpha,d_0}C_{\alpha,d_0}\nep^{\alpha\varphi(d_0)}+C_{\alpha,d_0}\Mn_{x\sim\pi}\left[\nep^{\psi(d(x,x_0))}\un_{d(x,x_0)\geq d_0}\right].$$
And then, since $C_{\alpha,d_0}<1$, we finally obtain:
$$\Mn_{x\sim\pi}\left[\nep^{\psi(d(x,x_0))}\un_{d(x,x_0)\geq d_0}\right]\leq \frac{C'_{\alpha,d_0}C_{\alpha,d_0}}{1-C_{\alpha,d_0}}\nep^{\alpha\varphi(d_0)}.$$
Now we just have to use the Markov inequality to derive the desired inequality.$\square$

\begin{rmk} In the previous proof, we didn't fully use the hypothesis $F(d_0)>\frac{s^2K(d_0)}{2}$.
In fact, for a fixed $d_0$, $\ln(C_{\alpha, d_0})$ is a convex function of $\alpha$ on the interval $[0,\frac{1}{s^2K(d_0)})$.
We have $C_{0,d_0}=1$ and $\frac{\partial}{\partial\alpha}\ln(C_{\alpha,d_0})|_{\alpha=0}<0$ if and only if $F(d_0)>\frac{s^2K(d_0)}{2}$.
So if $0<F(d_0)\leq\frac{s^2K(d_0)}{2}$, there doesn't exist any $\alpha$ such that $C_{\alpha,d_0}<1$ and so the theorem wouldn't tell us anything at all.
\end{rmk}

\begin{rmk} If $K(d_0)\leq\frac{1}{2}$ have $C_{\frac{2}{s^2},d_0}\geq1$, so we must have $\alpha<\frac{2}{s^2}$.
Under the hypothesis that $\kappa(x)\xrightarrow[x\rightarrow\infty]{}0$ and $F(x)\xrightarrow[x\rightarrow\infty]{}+\infty$, we can find for any $0<\alpha<\frac{2}{s^2}$ a $d_0$ such that $C_{\alpha,d_0}<1$.
Of course we need a greater $d_0$ when $\alpha$ gets closer to $\frac{2}{s^2}$.
\end{rmk}

One way to choose $\alpha$ and $d_0$ is given by the following proof of Theorem \ref{tm1}:

\tit{Proof of Theorem \ref{tm1}:} we use Theorem \ref{princ} with $\alpha=\frac{1}{2s^2}$ and $d_0=2\varepsilon+\frac{\ln(2)s^2}{\rho}$.
We only have to check that in this case, $C_{\alpha,d_0}\leq\nep^{-\frac{\rho^2}{4s^2}}$ and $C'_{\alpha,d_0}\leq\nep^{\frac{3\varepsilon}{2s^2}\max(3\varepsilon,\rho+\frac{\ln(2)s^2}{\rho})}$.

We have
$$\ln(C_{\alpha,d_0})=\left(-\alpha+\frac{\alpha^2s^2}{2(1-\alpha s^2K(d_0))}\right)F^2(d_0)-\frac{1}{2}\ln(1-\alpha s^2K(d_0)).$$
Since $K(d_0)\leq 1$, we have $-\alpha+\frac{\alpha^2s^2}{2(1-\alpha s^2K(d_0))}\leq-\frac{1}{4s^2}$.
We have $F(d_0)\geq\rho+\frac{\ln(2)s^2}{\rho}K(d_0)$, and then $F(d_0)^2\geq\rho^2+2\ln(2)s^2K(d_0)$.
Using the concavity of $\ln$ on $[\frac{1}{2},1]$, we get $\ln(1-\alpha s^2K(d_0))\geq-\ln(2)K(d_0)$.
Thus we get:
$$\ln(C_{\alpha,d_0})\leq-\frac{1}{4s^2}(\rho^2+2\ln(2)s^2K(d_0))+\frac{\ln(2)}{2}K(d_0)=-\frac{\rho^2}{4s^2}.$$

For $C'_{\alpha,d_0}$, we have 
\begin{align*}\ln(C'_{\alpha,d_0})&=\frac{1}{2s^2}\un_{J(x_0)+\varepsilon>d_0-F(d_0)}\int_{d_0-F(d_0)}^{J(x_0)+\varepsilon}\max(F(d_0),F(u))\der u\\
&\leq\frac{1}{2s^2}((J(x_0)+\varepsilon)-(d_0-F(d_0)))_+\max(F(d_0),F(J(x_0)+\varepsilon)).\end{align*}
By the triangular inequality for $W_1$, we have $J(x_0)\leq W_1(\delta_{x_0},P(x))+W_1(P(x),P(x_0))\leq W_1(\delta_{x_0},P(x))+d(x_0,x)$, for any $x$ because the coarse Ricci curvature is nonnegative.
If we take $x$ such that $\varepsilon<d(x_0,x)\leq2\varepsilon$, we have $J(x_0)\leq2d(x_0,x)-\rho\leq4\varepsilon-\rho$.
We have $F(d_0)\leq\rho+\frac{\ln(2)s^2}{\rho}$, so $d_0-F(d_0)\geq2\varepsilon-\rho$ and then $((J(x_0)+\varepsilon)-(d_0-F(d_0)))_+\leq3\varepsilon$.
And finally, $F(J(x_0)+\varepsilon)\leq F(5\varepsilon-\rho)\leq 3\varepsilon$.
Putting that together give us the desired bound for $C'_{\alpha,d_0}$.$\square$

\tit{Proof of Proposition \ref{tm2}:}
We take $\alpha$ and $d_0$ as in the proof of Theorem \ref{tm1}.
We consider the sequence of probability measures $P^n_{x_0}$.
Then, doing as in the proof of Theorem \ref{princ}, we have :
$$\Mn_{x\sim P^{n+1}_{x_0}}[\nep^{\psi(d(x,x_0))}\un_{d(x,0)\geq d_0}]\leq C'_{\alpha,d_0}C_{\alpha,d_0}+C_{\alpha,d_0}\Mn_{x\sim P^n_{x_0}}[\nep^{\psi(d(x,x_0))}\un_{d(x,0)\geq d_0}].$$
From that, we can conclude that there exists $C<+\infty$ such that for all $n$, we have $\Mn_{P^n_{x_0}}[\nep^{\psi(d(x,x_0))}]<C$.
So the sequence $P^n_{x_0}$ is tight, and then, so is the sequence $\pi_n=\frac{1}{n+1}\sum_{i=0}^nP^i_{x_0}$.
Because closed balls are compact, we can extract a weakly convergent subsequence $\pi_{\theta(n)}$, and we denote by $\pi$ its limit.
The $W_1$ distance metrizes the weak convergence on the set of probability measures on $X$ satisfying $\Mn[\nep^{\psi(d(x,x_0))}]<C$ (see \cite{VILL09}).
Thus the subsequence $\pi_{\theta(n)}$ converges to $\pi$ for the $W_1$ distance.
Furthermore, we have $W_1(\pi_n,P\pi_n)\leq\frac{C'}{n+1}$ with $C'<\infty$ a constant.
We have then
$$W_1(\pi,P_\pi)\leq W_1(\pi,\pi_{\theta(n)})+W_1(\pi_{\theta(n)},P\pi_{\theta(n)})+W_1(P\pi_{\theta(n)},P\pi).$$
The nonnegative coarse Ricci curvature implies that $P$ contracts the $W_1$ distance (\cite{OLLI09}), so the third term of the right hand side is at most the first one.
We have already seen that the first two terms tend to $0$ when $n$ tends to $+\infty$.
So the right hand side tends to $0$ when $n$ tends to $+\infty$.
Thus $W_1(\pi,P\pi)=0$, and then $\pi$ is an invariant measure.$\square$

\section{Some examples}

Let us see which concentration we can get with Theorem \ref{tm1} and Theorem \ref{princ} in some examples below.

\begin{xpl}[Discrete time M/M/k queue (see, for example \cite{GRST03})]
Let $0<n_0<k$ be two integers.
We consider here the Markov chain on integers with transition kernel:
\begin{align*} p(n,n+1)&=\frac{n_0}{n_0+k}\\
p(n,n)&=\frac{(k-n)_+}{n_0+k}\\
p(n,n-1)&=\frac{\min(n,k)}{n_0+k}\\
p(n,m)&=0\quad\textrm{if }|n-m|>1.\\
\end{align*}
\end{xpl}
The origin $x_0$ we will consider to apply Theorem \ref{tm1} is $n_0$, the only point at which the probability to jump at left equals the probability to jump at right (that is why we chose $n_0$ integer).
Hoeffding's Lemma (see \cite{HOEF63}) states that for a random variable $X$ such that $a\leq X\leq b$ almost surely, we have $\Mn\left[\nep^\lambda(X-\Mn[X])\right]\leq\nep^\frac{\lambda^2(b-a)^2}{8}$.
So we can take $s=1$ in theorem \ref{tm1}.
To compute the coarse Ricci curvature, we remark that if $x<y$, the measure $P_y$ dominates stochastically the measure $P_x$, and thus the $W_1$ distance between them is the difference of their expectations.
For $x<y$, the coarse Ricci curvature $K(x,y)$ is then $\frac{1}{n_0+k}$ if $y\leq k$, $\frac{k-x}{y-x}\frac{1}{n_0+k}$ if $x<k<y$ and $0$ if $x\geq k$.
If we take $\varepsilon=1$, we have $\rho=\frac{1}{n_0+k}$, and $K(r)=\frac{\un_{r<k-n_0}}{n_0+k}$.

Applying Theorem \ref{tm1} should give a Gaussian then exponential concentration, but, as $\rho$ is very small, $d_0$ is large $(2+(n_0+k)\ln(2))$.
If $k-n_0\leq\frac{2\ln(2)n_0+2}{1-\ln(2)}$, we get only the exponential part.
If $k$ is too large, $d_0$ is large too, and the gaussian-then-exponential bounds starts far away from $n_0$.
We can try to take a larger $\varepsilon$ to get a better $\rho$.
Indeed, we get $\rho=\frac{\min(\varepsilon,k-n_0)}{n_0+k}$, but we pay that by a worse curvature $K(r)=\frac{1}{n_0+k}\min(1,\max(0,\frac{k-n_0-r}{\varepsilon}))$.
We distinguish 3 cases depending on how $k-n_0$ is tall with respect to $n_0$.

When $k-n_0$ is between $\sqrt{n_0}$ and $n_0$, the equilibrium measure is well approximated by a Gaussian between $0$ and $k$.
\begin{center}
\includegraphics{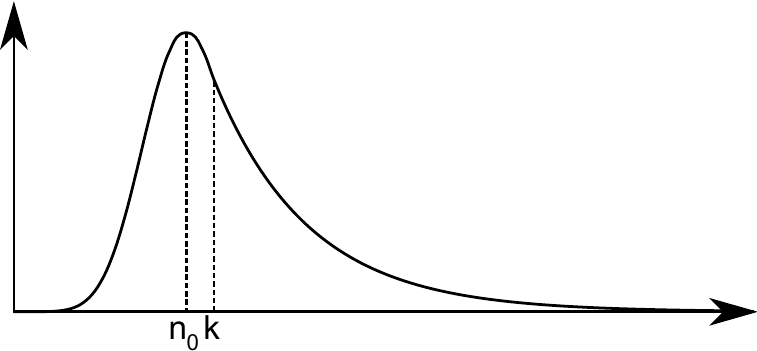}
\end{center}
The optimal $\varepsilon$ is $O(\sqrt{n_0})$, the coefficient of the Gaussian part of the concentration inequality is $O(\frac{1}{n_0})$, which is good, and the coefficient of the exponential part is $O(\frac{k-n_0}{n_0})$, like the right one.

When $k-n_0$ is $o(\sqrt{n_0})$, the mass of $[0,k]$ under the equilibrium measure is negligible with respect to the one of $[k,\infty)$.
\begin{center}
\includegraphics{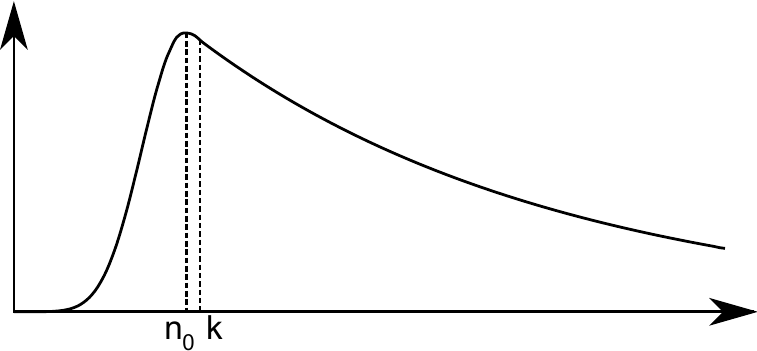}
\end{center}
The optimal $\varepsilon$ and $d_0$ are $O(k-n_0)$, this time, we have no Gaussian part because $d_0$ is too large (and indeed, there is no Gaussian part in the equilibrium measure), and the coefficient of the exponential part is about one half of the right one.

When $k-n_0$ is greater than $n_0$, the equilibrium measure is almost the Poissonian one with parameter $n_0$, the density of the equilibrium measure is illustrated below:
\begin{center}
\includegraphics{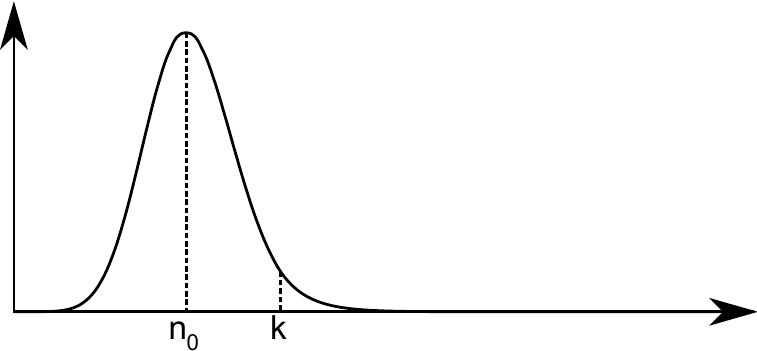}
\end{center}
The optimal $\varepsilon$ and $d_0$ are $O(\sqrt{k})$, the coefficient appearing in the Gaussian part is $O(\frac{1}{k})$, instead of an expected $\frac{1}{n_0}$, and the coefficient of the exponential part is $O(1)$, which is clearly not optimal, so Theorem \ref{tm1} gives a rather bad concentration inequality.

\begin{xpl}[Discrete time Ornstein Uhlenbeck]
Let $0<\alpha\leq1$ be a real parameter.
Here we consider the Markov Chain on $\Rl$ given by the transition kernel:
$$P_x=\mathcal{N}((1-\alpha)x,1).$$
\end{xpl}

It is shown in \cite{BOBGT99} that in the Gaussian case, we can take the variance of the distribution for $s^2$.
So we take $s^2=1$, and for every $\varepsilon>0$, the curvature is constant $K=\alpha$.
We have $\rho=\varepsilon-\sqrt{\frac{2}{\pi}}\left(\nep^{-\frac{(1-\alpha)^2\varepsilon^2}{2}}+\int_0^{(1-\alpha)\varepsilon}\nep^{-\frac{x^2}{2}}\der x\right)\geq-\sqrt{2}{\pi}+\alpha\varepsilon$.
Theorem \ref{tm1} applied with $\varepsilon=\frac{\sqrt{2\ln(2)\alpha}+\sqrt{{\frac{8}{\pi}}}}{\alpha}$ gives us Gaussian concentration with coefficient $\frac{\alpha}{4}$ instead of $\frac{\alpha(2-\alpha)}{2}$ (so we have a loss of a factor between $2$ and $4$), and $d_0=O(\sqrt{\frac{1}{\alpha}})$

The bad behaviour of $s^2$ prevents to easily generalize Theorem \ref{tm1} or Theorem \ref{princ} to continuous time.
The following example of a continuous time processes, whose generator merges a diffusive part and a jump part, shows that a generalization of those theorems does not hold, even if the jump rate is uniformly bounded.

\begin{xpl}
Consider a continuous time process on $\Rl_+$ with a linear drift towards $0$ and a random jump to the right of size $1$ and rate $1$.
The generator of this process is given by $Lf(x)=-\alpha x f'(x)+f(x+1)-f(x)$, with $\alpha>0$ a constant which quantifies the drift.
\end{xpl}

In this example, the coarse Ricci curvature is $\alpha$.
Indeed, using the coupling of the processes $X_t$ and $Y_t$ starting at $x$ and $y$ such that $X_t$ and $Y_t$ jump at the same times shows that the law of $Y_t$ is the translation of the law of $X_t$ by $(y-x)\nep^{-\alpha t}$.
If something like Theorem \ref{tm1} or Theorem \ref{princ} did hold, we would have Gaussian concentration.
But actually there is only Poissonian concentration.
Let us prove there is Poissonian concentration and no better.
We denote by $X_t$ the value of the process at the time $T$.
Let $T_1,T_2,\dots$ be the successive times of the jumps.
For all $T>0$, let $N(T)$ be the number of jumps between $0$ and $T$.
We have
$$X_T=\nep^{-\alpha T}X_0+\sum_{i=1}^{N(T)}\nep^{-\alpha(T-T_i)}.$$
If we take $X_0=0$ then $\Mn[X_T]\leq T$, and since the coarse Ricci curvature is greater than $\alpha>0$, there exists an unique invariant probability measure $\pi$ (see \cite{VEY12}).

Now take $X_0$ with the law $\pi$.
Then $X_1$ has the law $\pi$, and is greater than $\nep^{-\alpha}N(1)$, which has a Poissonian concentration since $N(1)$ follows precisely a Poisson law of parameter $1$.
So we cannot have a better concentration than a Poissonian one.

It remains to prove that $\pi$ has Poissonian concentration.
We take $X_0=0$.
Let us consider the Laplace/Fourier transform of $X_T$, that is $G_T(\lambda)\coleq\Mn[\nep^{\lambda X_T}]$ for $\lambda\in\C$.
$N(T)$ has the law $\mathcal{P}(T)$, and the repartition of the $T_i$'s knowing $N(T)$ is the one of $N(T)$ independent random variables uniformly distributed in $[0,T]$.
So we have:
\begin{align*}G_T(\lambda)&=\sum_{k=0}^\infty\frac{T^k\nep^{-T}}{k!}\left(\frac{1}{T}\int_0^T\nep^{\lambda\nep^{-\alpha t}}\der t\right)^k\\
&=\sum_{k=0}^\infty\frac{\nep^{-T}}{k!}\left(\int_0^T\left(\sum_{n=0}^\infty\frac{\lambda^n\nep^{-n\alpha t}}{n!}\right)\der t\right)^k\\
&=\nep^{\sum_{n=1}^\infty\left[\frac{-\lambda^n\nep^{-n\alpha t}}{n\alpha n!}\right]_{t=0}^{T}}\\
&=\nep^{\frac{I(\lambda)-I(\lambda\nep^{-\alpha T})}{\alpha}}\end{align*}
with $I(\lambda)=\sum_{n=1}^\infty\frac{\lambda^n}{nn!}=\int_0^\lambda\frac{(\nep^z-1)\der z}{z}$.
We see that $G_T(\lambda)$ tends to a limit $G(\lambda)$, which is the Laplace/Fourier transform of $\pi$, when $T$ tends to $+\infty$.

We have $G(\lambda)=\nep^{\frac{I(\lambda)}{\alpha}}$.
An integration by parts gives us $I(\lambda)=\frac{\nep^{\lambda}-\lambda-1}{\lambda}+\int_0^\lambda\frac{\nep^z-z-1}{z^2}\der z$, so $I(\lambda)\sim\frac{\nep^\lambda}{\lambda}$.
For $l>1$, we use the Markov inequality on $\nep^{\ln(l)X}$ and get:
$$\Pb_\pi[X\geq l]\leq\nep^{\frac{I(\ln(l))}{\alpha}-l\ln(l)},$$
and we have $\frac{I(\ln(l))}{\alpha}\sim\frac{l}{\alpha\ln(l)}=o(l\ln(l))$.
So we have Poissonian concentration.
\bibliographystyle{abbrv}
\bibliography{base}
\end{document}